\documentclass[10pt]{article}
\usepackage[francais,english]{babel}
\usepackage[latin1]{inputenc}
\usepackage[T1]{fontenc}
\usepackage{amsmath}
\usepackage{amsfonts}
\usepackage{amssymb}
\usepackage{amsthm}
\newtheorem{prop}{Proposition}
\newtheorem{lem}[prop]{Lemme}
\newtheorem{sublem}[prop]{Sous-Lemme}
\newtheorem{rmk}[prop]{Remarque}
\newcommand{\Spec}{\mathop\mathrm{Spec}\nolimits}
\begin{document}
\selectlanguage{francais}
\title{Les hypersurfaces cubiques sont\\séparablement rationnellement connexes}
\author{David A. Madore}
\maketitle
\begin{center}
CVS: \verb=$Id: freecurve.tex,v 1.11 2006-05-25 14:03:16 david Exp $=
\end{center}

\selectlanguage{english}
\begin{abstract}
This note (which makes no claim to novelty) presents a proof of the
separable rational connectedness of smooth cubic hypersurfaces, in any
characteristic, by showing how to explicitly construct very free
curves (of degree $3$) on them.
\end{abstract}
\selectlanguage{francais}
\begin{abstract}
Cette note (qui ne prétend pas à l'originalité) démontre la séparable
rationnelle connexité des hypersurfaces cubiques lisses, en toute
caractéristique, en construisant explicitement des courbes très libres
(de degré $3$) tracées dessus.
\end{abstract}

\par\bigskip

\par\noindent\textbf{Introduction :} L'objet de cette note est de
fournir une démonstration simple d'un fait qui doit certainement être
considéré comme connu mais qui semble difficile à trouver dans la
littérature : le fait que, \emph{en toute caractéristique}, les
hypersurfaces cubiques lisses sont \emph{séparablement} rationnellement
connexes.  On envoie à \cite{KollarBook} (notamment chap. IV) ou
\cite{Debarre} pour une discussion générale sur la (séparable)
rationnelle connexité ; ici on prendra pour définition l'existence
d'une courbe \emph{très libre} : c'est-à-dire qu'on dira qu'une
variété projective lisse (intègre) $X$ sur un corps $k$ algébriquement
clos est séparablement rationnellement connexe lorsqu'il existe
$h\colon \mathbb{P}^1\to X$ tel que le fibré $h^* T_X$ soit ample
(i.e., $H^1(\mathbb{P}^1,(h^* T_X)(-2))=0$), où $T_X$ désigne le fibré
tangent à $X$.  Pour les considérations générales sur les
hypersurfaces cubiques, et notamment la rationalité des surfaces
cubiques, on renvoie à \cite{Manin}.

La question plus générale de la séparable rationnelle connexité, en
toute caractéristique, des hypersurfaces de Fano lisses, c'est-à-dire
des hypersurfaces lisses de degré $d$ dans $\mathbb{P}^n$ avec $n\geq
d$, semble encore ouverte (comparer \cite{KollarBook} V.2.13
et V.5.11).

\bigskip

\begin{prop}
Soit $X$ une variété projective lisse (intègre) sur $k$ un corps
algébriquement clos, et soit $D$ une sous-variété lisse (intègre) de
codimension $1$ qui, vue en tant que diviseur sur $X$ (ayant une seule
composante, avec multiplicité $1$), est ample.  On se donne $h \colon
\mathbb{P}^1 \to D$ non constant, et on suppose que $h$ est très libre
à valeurs dans $D$, c'est-à-dire $H^1 (\mathbb{P}^1, T_D(-2)) = 0$.
Alors $h \colon \mathbb{P}^1 \to X$ est encore très libre, vue comme
courbe sur $X$, c'est-à-dire $H^1 (\mathbb{P}^1, T_X(-2)) = 0$.
\label{FreenessPropOne}
\end{prop}
\begin{proof}
Soit $i\colon D\to X$ le morphisme d'immersion de $D$ dans $X$.  On a
la suite exacte
$$
0 \to T_D \to i^* T_X \to \mathcal{O}_D(D) \to 0
$$
où $T_D$ est le fibré tangent à $D$, $T_X$ celui à $X$,
$\mathcal{O}_X(D)$ le fibré en droites (supposé ample) associé au
diviseur $D$, et $\mathcal{O}_D(D)$ sa restriction à $D$ lui-même, qui
n'est autre que le fibré normal à $D$.  En tirant cette suite exacte
courte de fibrés par le morphisme $h \colon \mathbb{P}^1 \to D$, on
trouve :
$$
0 \to h^* T_D \to h^* T_X \to \mathcal{O}_{\mathbb{P}^1}(\deg_D h) \to 0
$$
où $\deg_D h$ désigne le degré de $h$ mesuré par rapport au diviseur
$D$ : on a $\deg_D h > 0$ car $D$ est ample et que $h$ est non
constant.  On en déduit en particulier, au niveau de la cohomologie,
pour tout $\ell \in\mathbb{Z}$ :
$$
H^1 (\mathbb{P}^1, (h^* T_D)(\ell))
\to H^1 (\mathbb{P}^1, (h^* T_X)(\ell))
\to H^1 (\mathbb{P}^1, (\mathcal{O}_{\mathbb{P}^1}(\deg_D h + \ell))) \to 0
$$
Prenons $\ell = -2$ : alors $H^1 (\mathbb{P}^1, (h^* T_D)(-2))$
s'annule par hypothèse, et $H^1 (\mathbb{P}^1, \penalty0
(\mathcal{O}_{\mathbb{P}^1}(\deg_D h - 2)))$ s'annule car $\deg_D h-2
\geq -1$.  Il s'ensuit qu'on a $H^1 (\mathbb{P}^1,\penalty0 (h^*
T_X)(-2)) = 0$, ce qu'on voulait démontrer.
\end{proof}

\begin{prop}
Soit $X$ une surface projective lisse (intègre) sur $k$ un corps
algébriquement clos, et $C\subseteq X$ une courbe rationnelle intègre
dans $X$, n'ayant pas d'autre singularité que des points doubles
ordinaires.  On suppose que $\deg_{-K_X} C \penalty0 \geq 3$, où
$\deg_{-K_X} C$ désigne le degré d'intersection de $C$ par rapport au
diviseur anticanonique $-K_X$ sur $X$.  Soit $h \colon \mathbb{P}^1
\to C$ la normalisation.  Alors $h \colon \mathbb{P}^1 \to X$ est très
libre, vue comme courbe sur $X$, c'est-à-dire $H^1 (\mathbb{P}^1,
T_X(-2)) = 0$.  (Réciproquement, si $h$ est très libre, alors
$\deg_{-K_X} C \penalty0 \geq 3$.)
\label{FreenessPropTwo}
\end{prop}
\begin{proof}
Rappelons que $T_{\mathbb{P}^1} \cong \mathcal{O}_{\mathbb{P}^1}(2)$.
On a une flèche $T_{\mathbb{P}^1} \to h^* T_X$ déduite par dualité de
$h^* \Omega^1_{X/\Spec k} \to \Omega^1_{\mathbb{P}^1_k/\Spec k}$, et
cette flèche est injective sur chaque fibre car la différentielle de
$h$ ne peut pas s'annuler (en un point lisse, c'est évident, et
au-dessus d'un point double, $h$ prend la direction d'une des deux
tangentes distinctes en ce point).  La flèche $I_C/I_C^2 \to i^*
\Omega^1_{X/\Spec k}$ (où $I_C$ est le faisceau d'idéaux définissant
$C$ sur $X$ et $i$ la flèche d'immersion de $C$ dans $X$) donne, en
tirant par $h$ et en dualisant, une flèche $h^* T_X \to
\mathcal{O}_{\mathbb{P}^1}(\deg_C h)$ (non nécessairement
surjective !).  La composée $T_{\mathbb{P}^1} \to h^* T_X \to
\mathcal{O}_{\mathbb{P}^1}(\deg_C h)$ s'annule, car elle s'annule sur
chaque fibre au-dessus de tout l'ouvert de lissité de $C$.  Mais en
regardant fibre à fibre, on voit même que la suite $T_{\mathbb{P}^1}
\to h^* T_X \to \mathcal{O}_{\mathbb{P}^1}(\deg_C h)$ est exacte : en
effet, le noyau de la flèche de droite est un sous-fibré de rang $1$
dans $h^* T_X$ qui contient l'image de la flèche de gauche, image qui,
comme on l'a expliqué, a rang $1$ en chaque fibre, donc il y a
exactitude fibre à fibre, d'où la suite exacte :
$$
0 \to \mathcal{O}_{\mathbb{P}^1}(2)
\to h^* T_X \to \mathcal{O}_{\mathbb{P}^1}(\deg_C h)
$$
La dernière flèche de cette suite exacte n'est pas nécessairement
surjective (elle ne l'est d'ailleurs jamais si $C$ a effectivement des
singularités), mais l'image de cette flèche est un sous-fibré de
$\mathcal{O}_{\mathbb{P}^1}(\deg_C h)$, manifestement de rang $1$ et
de degré celui de $h^* T_X$ moins deux (car le degré de $h^* T_X$ est
$\deg_{-K_X} C$), c'est-à-dire $\deg_{-K_X} C - 2$, qui est donc au
moins égal à $1$ par hypothèse.  On peut donc écrire
$$
0 \to \mathcal{O}_{\mathbb{P}^1}(2)
\to h^* T_X \to \mathcal{O}_{\mathbb{P}^1}(\deg_{-K_X} C - 2) \to 0
$$
(notons au passage que si $C$ est lisse, on sait que $\deg_C h =
\deg_{-K_X} C - 2$ puisque le genre de $C$ est zéro, donc la dernière
flèche de la première suite exacte était bien surjective ; et
réciproquement).  En passant comme dans la démonstration précédente à
la suite exacte longue de cohomologie (ce sont les $H^1$ qui nous
intéressent) après avoir tensorisé par
$\mathcal{O}_{\mathbb{P}^1}(-2)$, on voit que $H^1 (\mathbb{P}^1,
T_X(-2)) = 0$ comme on le voulait.

Pour la réciproque, on écrit $h^* T_X \cong \mathcal{O}_{\mathbb{P}^1}
(d_1) \oplus \mathcal{O}_{\mathbb{P}^1} (d_2)$ avec $d_1 \geq d_2$ :
l'hypothèse que $h$ est très libre se traduit $d_2\geq 1$, et
l'existence d'une flèche injective $\mathcal{O}_{\mathbb{P}^1}(2) \to
h^* T_X$ donne $d_2\geq 2$, donc $\deg_{-K_X} C = d_1 + d_2 \geq 3$,
la conclusion recherchée.
\end{proof}

En particulier, toute cubique plane intègre à point double ordinaire
tracée sur une surface cubique lisse est très libre :
\begin{prop}
Soit $X \subseteq \mathbb{P}^3$ une surface cubique lisse sur un corps
$k$ algébriquement clos (de caractéristique arbitraire), $C \subseteq
X$ la courbe cubique plane intersection de $X$ avec un plan $\Pi$ tel
que $C$ ait un point singulier $x$ (c'est-à-dire que $\Pi$ est le plan
$\Pi(x)$ tangent à $X$ en $x$) et supposons que $C$ soit intègre et
que $x$ soit un point double ordinaire (par opposition à un cusp) :
alors la normalisation $h \colon \mathbb{P}^1 \to C$ composée avec
l'inclusion canonique définit une courbe $h \colon \mathbb{P}^1 \to X$
très libre sur $X$.
\label{CubSurfProp}
\end{prop}
\begin{proof}
On se trouve dans les conditions d'application de la
proposition \ref{FreenessPropTwo} : le fibré anticanonique $-K_X$ sur
une surface cubique est donné par une section plane ; le degré de la
courbe $C$ d'intersection est alors $3$, ce qui assure que les
hypothèses sont bien vérifiées.
\end{proof}
\begin{proof}[Démonstration par calcul explicite.]
On appelle $(X_0:X_1:X_2:X_3)$ les coordonnées de $\mathbb{P}^3$, et
on suppose (sans perte de généralité) que $X_3 = 0$ est l'équation du
plan $\Pi$.  La courbe $C$ est une cubique rationnelle, donc elle a un
(unique) point singulier $x$, qui, par hypothèse, est un point double
ordinaire.  Mettons que ce point soit donné par $X_1 = X_2 = 0$ (soit
$(1:0:0)$) dans le plan $\Pi$.  L'équation de $C$ s'écrit alors $X_0\,
q(X_1, X_2) + c(X_1, X_2) = 0$ où $q$ est une forme quadratique en
$X_1, X_2$ et $c$ une forme cubique.  L'hypothèse que $C$ a en
$(1:0:0)$ un point double ordinaire se traduit par le fait que $q$ est
non nulle et a deux racines distinctes (correspondant aux deux
directions tangentes), et l'hypothèse que $C$ est intègre, donc ne
contient aucune droite, se traduit par le fait que $q$ et $c$ sont
sans racine commune.  Quitte à faire un changement de coordonnées, on
peut supposer que $q(X_1, X_2) = X_1 X_2$.  Écrivons $c(X_1,X_2) =
\alpha_0 X_1^3 + \alpha_1 X_1^2 X_2 + \alpha_2 X_1 X_2^2 + \alpha_3
X_2^3$, où $\alpha_0\neq 0$ et $\alpha_3\neq 0$ (sans quoi $q$ et $c$
auraient un zéro commun).  Quitte à remplacer $X_0$ par $X_0 +
\alpha_1 X_1 + \alpha_2 X_2$, on peut supposer que $\alpha_2=0$ et
$\alpha_3=0$ ; et quitte à multiplier $X_1$ par $\sqrt{\alpha_0}$ et
$X_2$ par $\sqrt{\alpha_3}$, on peut de plus prendre $\alpha_0=1$ et
$\alpha_3=1$.  La normalisation $h \colon \mathbb{P}^1 \to C$ est
alors donnée par le paramétrage $(U:V) \mapsto (-U^3-V^3 : U^2 V : U
V^2)$.  On peut alors écrire $h^* T_{\mathbb{P}^3} \cong
\mathcal{O}_{\mathbb{P}^1}(5) \xi \oplus \mathcal{O}_{\mathbb{P}^1}(4)
\eta \oplus \mathcal{O}_{\mathbb{P}^1}(3) \frac{\partial}{\partial
X_3}$ avec $\xi \in \Gamma (\mathbb{P}^1, (h^* T_{\mathbb{P}^3})
(-5))$ donné par $\xi = \frac{U^2}{V^4} \frac{\partial}{\partial X_0}
- \frac{U}{V^3} \frac{\partial}{\partial X_1} - \frac{1}{V^2}
\frac{\partial}{\partial X_2} = \frac{V^2}{U^4}
\frac{\partial}{\partial X_0} - \frac{1}{U^3} \frac{\partial}{\partial
X_1} - \frac{V}{U^3} \frac{\partial}{\partial X_2}$ et $\eta \in
\Gamma (\mathbb{P}^1, (h^* T_{\mathbb{P}^3})\penalty0 (-4))$ donné par
$\eta = \frac{U}{V^2} \frac{\partial}{\partial X_0} -
\frac{1}{V}\frac{\partial}{\partial X_1} = - \frac{V}{U^2}
\frac{\partial}{\partial X_0} + \frac{1}{U}\frac{\partial}{\partial
X_2}$.

Soit $f = X_0 X_1 X_2 + c(X_1, X_2) + X_3 Q(X_0,X_1,X_2) + X_3^2
L(X_0,X_1,X_2) + A X_3^3$ l'équation de la surface $X$, où $Q$ est une
forme quadratique telle que $Q(1,0,0) \neq 0$ (sans quoi $X$ serait
singulière en $(1:0:0:0)$), $L$ une forme linéaire (éventuellement
nulle), et $A$ une constante.  On peut alors calculer $\xi\cdot f =
-U^2 V^2$ et $\eta\cdot f = - U^4 V - U V^4$ et
$\frac{\partial}{\partial X_3} f = Q(-U^3-V^3, U^2 V, U V^2)$.  Ainsi,
aucune combinaison linéaire des sections (de $\Gamma (\mathbb{P}^1,
(h^* T_{\mathbb{P}^3})\penalty0 (-3))$) $U^2\xi$, $UV\xi$, $V^2\xi$,
$U\eta$, $V\eta$ ou $\frac{\partial}{\partial X_3}$ n'annule $f$,
c'est-à-dire, n'est dans le noyau $h^* T_X$ de la flèche $h^*
T_{\mathbb{P}^3} \to \mathcal{O}_{\mathbb{P}^1} (9)$ : en effet, en
faisant $V=0$ on se convainc que le coefficient devant
$\frac{\partial}{\partial X_3}$ doit être nul, et l'annulation des
autres coefficients est claire.  On a ainsi prouvé $\Gamma
(\mathbb{P}^1, (h^* T_X)\penalty0 (-3)) = 0$ ; comme $h^* T_X$ doit
s'écrire $\mathcal{O}_{\mathbb{P}^1} (d_1) \oplus
\mathcal{O}_{\mathbb{P}^1} (d_2)$ où mettons $d_1\geq d_2$, ceci
signifie $d_1 < 3$, mais comme $d_1+d_2 = 3$, on a manifestement
$d_1=2$ et $d_2=1$.  Ceci montre bien $H^1 (\mathbb{P}^1, (h^* T_X)
(-2)) = 0$.  (On peut aussi dire que le
$\mathcal{O}_{\mathbb{P}^1}(2)$ de la décomposition de $h^* T_X$ est
engendré par $UV\eta - U^3\xi + V^3 \xi = - 3 \frac{V^2}{U}
\frac{\partial}{\partial X_0} + U \frac{\partial}{\partial X_1} + 2V
\frac{\partial}{\partial X_2} = 3 \frac{U^2}{V}
\frac{\partial}{\partial X_0} - 2U \frac{\partial}{\partial X_1} -V
\frac{\partial}{\partial X_2}$, qui, de fait, correspond bien à
l'image de la flèche $\mathcal{O}_{\mathbb{P}^1}(2) \to h^* T_{X}$
explicitée dans la démonstration précédente.)
\end{proof}

À l'opposé, une cubique plane intègre \emph{cuspidale} tracée sur une
surface cubique lisse \emph{n'est pas} très libre (et l'hypothèse
faite dans la proposition \ref{FreenessPropTwo} sur les singularités
de $C$ est donc essentielle) :

\begin{rmk}
Soit $X \subseteq \mathbb{P}^3$ une surface cubique lisse sur un corps
$k$ algébriquement clos (de caractéristique arbitraire), $C \subseteq
X$ la courbe cubique plane intersection de $X$ avec un plan $\Pi$ tel
que $C$ ait un point singulier $x$ (c'est-à-dire que $\Pi$ est le plan
$\Pi(x)$ tangent à $X$ en $x$) et supposons que $C$ soit intègre et
que $x$ soit un cusp : alors la normalisation $h \colon \mathbb{P}^1
\to C$ composée avec l'inclusion canonique définit une courbe $h
\colon \mathbb{P}^1 \to X$ qui n'est pas très libre sur $X$.
\end{rmk}
\begin{proof}
On reprend les notations utilisées dans le calcul explicite prouvant
la proposition \ref{CubSurfProp}, mais cette fois avec $q(X_1,X_2) =
X_2^2$.  Si la caractéristique est différente de $3$, quitte à faire
des changements linéaires de coordonnées, on peut supposer que
$c(X_1,X_2) = X_1^3$.  Alors l'existence de la section $\delta \in
\Gamma (\mathbb{P}^1, (h^* T_X) (-3))$ donnée par $\delta = 3
\frac{U^2}{V^2} \frac{\partial}{\partial X_0} -
\frac{\partial}{\partial X_1} = 2 \frac{\partial}{\partial X_1} + 3
\frac{V}{U} \frac{\partial}{\partial X_2}$ montre que $d_1\geq 3$ (en
fait, on a précisément $h^* T_X \cong \mathcal{O}_{\mathbb{P}^1} (3)
\oplus \mathcal{O}_{\mathbb{P}^1}$, ce qui se vérifie par l'annulation
de $\Gamma (\mathbb{P}^1, (h^* T_X) (-4))$, mais ce n'est pas
nécessaire pour savoir que la courbe n'est pas très libre).  En
caractéristique $3$, on écrit $c(X_1,X_2) = X_1^3 + \alpha X_1^2 X_2$,
et l'expression de la section $\delta \in \Gamma (\mathbb{P}^1, (h^*
T_X) (-3))$ est plus compliquée : $\delta = - \alpha \frac{U}{V}
\frac{\partial}{\partial X_0} - \frac{\partial}{\partial X_1} =
\alpha^2 \left(\alpha\frac{V}{U} - 1\right) \frac{\partial}{\partial
X_0} - \left(\alpha\frac{V}{U} - 1\right)^2 \frac{\partial}{\partial
X_1} - \alpha \frac{V^2}{U^2} \left(\alpha\frac{V}{U} + 1\right)
\frac{\partial}{\partial X_2}$, mais la conclusion est la même.
\end{proof}

{\footnotesize

\textbf{Note :} En analysant plus précisément la structure de la
démonstration que nous venons de faire, on voit qu'en fait la surface
$X$ dans laquelle $C$ est plongée importe peu : il y aura
\emph{toujours} une flèche non nulle $\mathcal{O}_{\mathbb{P}^1} (3)
\to h^* T_X$ donnée par $\delta$ ; et même, la fin de la démonstration
de la proposition \ref{FreenessPropTwo} tient encore si on remplace
essentiellement $\mathcal{O}_{\mathbb{P}^1} (2)$ par
$\mathcal{O}_{\mathbb{P}^1} (3)$ tout du long.  Ceci signifie que, si
$C$ (la cubique cuspidale plane, vue comme une courbe abstraite) est
plongée dans une surface projective lisse $X$, la flèche de
normalisation la présente comme très libre sur $X$ si et seulement si
$\deg_{-K_X} C \penalty0 \geq 4$.  L'auteur de cette note soupçonne qu'il
existe un invariant $\rho$ facilement calculable, pour une courbe
singulière $C$ (point singulier par point singulier, additif, et
valant $0$ pour un point double ordinaire et $1$ pour un cusp cubique
simple) tel que $C$ plongée dans une surface projective lisse $X$ soit
très libre si et seulement si $\deg_{-K_X} C \geq 3+\rho$.

}

\medbreak

Il n'y a bien sûr aucune surprise à l'existence de courbes très libres
sur une surface cubique lisse, puisque cette dernière est rationnelle,
donc certainement séparablement rationnellement connexe.  Le fait
intéressant dans la proposition \ref{CubSurfProp} est qu'on peut en
trouver explicitement.  Ces courbes existent d'ailleurs, en vertu du
lemme suivant (qui peut être intéressant en lui-même) :
\begin{lem}
Soit $X \subseteq \mathbb{P}^3$ une surface cubique lisse sur un corps
$k$ algébriquement clos de caractéristique différente de $2$.  Alors il
existe $x \in X$ tel que l'intersection $C(x)$ de $X$ avec le plan
tangent $\Pi (x)$ à $X$ en $x$ soit une courbe cubique intègre ayant
pour unique singularité un point double ordinaire en $x$.
\label{CubSurfLemma}
\end{lem}
\begin{sublem}
Soient $P_0,P_1,P_2,P_3,P_4,P_5$ six points (distincts) du plan
projectif sur un corps $k$ de caractéristique différente de $2$ : on
suppose que trois quelconques d'entre eux ne sont pas alignés et que
tous les six ne sont pas situés sur une même conique.  Alors il existe
un point $Q$ à l'intersection de deux des (quinze) droites reliant
deux des six points $P_i$ qui est distinct de tous les $P_i$, non
situé sur aucune autre des droites reliant deux des $P_i$ et également
non situé sur une des six coniques reliant cinq des $P_i$.
\label{SixPointSublemma}
\end{sublem}
\begin{proof}[Démonstration du sous-lemme]
Les points $P_0,P_1,P_2,P_3$ forment une base projective du plan.
Quitte à effectuer une transformation projective (qui ne change rien à
la situation ni à la conclusion recherchée), on peut donc supposer
qu'ils ont les coordonnées respectives $(0:0:1) ,\penalty-100 (1:0:1)
,\penalty-100 (1:1:1) ,\penalty-100 (0:1:1)$.  À ce moment-là, le
point d'intersection de $(P_0 P_1)$ et de $(P_2 P_3)$ est $(1:0:0)$,
celui de $(P_0 P_2)$ et de $(P_1 P_3)$ est $(1:1:2)$, et celui de
$(P_0 P_3)$ et de $(P_1 P_2)$ est $(0:1:0)$ : la chose qui nous
importe est qu'ils ne sont pas alignés en caractéristique différente
de $2$, donc ils ne peuvent pas être tous les trois sur la droite
$(P_4 P_5)$.  Appelons $Q$ un de ces trois points qui n'est pas situé
sur $(P_4 P_5)$, et mettons pour fixer les idées (et sans perte de
généralité) que ce soit le point d'intersection de $(P_0 P_1)$ et de
$(P_2 P_3)$.  Manifestement, $Q$ est distinct des six points $P_i$.
Il n'est situé sur aucune des droites $(P_0 P_2)$, $(P_0 P_3)$, $(P_0
P_4)$ ou $(P_0 P_5)$, sans quoi $P_0$ et $P_1$ seraient alignés avec
un $P_i$ pour $i\geq 2$.  Pour des raisons semblables, $Q$ ne peut pas
être situé sur $(P_1 P_2)$, $(P_1 P_3)$, $(P_1 P_4)$, $(P_1 P_5)$,
$(P_2 P_4)$, $(P_2 P_5)$, $(P_3 P_4)$ ou $(P_3 P_5)$, et on a déjà
expliqué qu'il n'était pas sur $(P_4 P_5)$.  Restent enfin les
coniques ; une conique définie par cinq des six points $P_i$ est
non-dégénérée (i.e., lisse), donc ne peut pas contenir trois points
alignés, et comme elle contient soit $P_0$ et $P_1$ soit $P_2$ et
$P_3$ elle ne peut pas contenir aussi $Q$.  Ce qui termine la
démonstration.
\end{proof}
\begin{proof}[Démonstration du lemme]
Comme $X$ est isomorphe à l'éclaté du plan projectif en six points, le
sous-lemme nous permet de trouver un point $z$ sur $X$ situé sur deux,
mais pas trois (on dit que $z$ n'est \emph{pas} un « point
d'Eckardt »), des vingt-sept droites tracées sur $X$.  (On peut aussi
faire appel au fait « connu » que la surface cubique ayant le plus de
points d'Eckardt est la surface de Clebsch, qui en a $10$, ce qui est
strictement moins que le tiers du nombre de points d'intersection,
soit $105$, de deux des $27$ droites tracées sur $X$.)  Ceci signifie
que $C(z)$ est formé de trois droites dans le plan $\Pi(z)$, se
croisant en trois points distincts.

Soit $D$ une des deux doites qui passent par $z$.  Considérons un
point $y$ sur $D$ : alors $C(y)$ est formé de la réunion de $D$ et
d'une conique (éventuellement dégénérée).  Comme il n'y a qu'un nombre
fini de droites tracées sur $X$, il y a un ouvert \emph{non vide} (de
$D$) de points $y$ pour lesquels la conique en question est lisse.
Par ailleurs, comme c'est le cas pour $z$, il y a un ouvert \emph{non
vide} de points $y$ pour lesquels elle rencontre la droite $D$ en deux
points distincts.  Fixons $y$ dans l'intersection de ces deux
ouverts : pour résumer, $C(y)$ est donc la réunion de $D$ et d'une
conique ${\mit\Gamma}$ qui coupe $D$ en $y$ et en un deuxième point
$y'$ distinct de $y$.

Considérons maintenant un point $x$ sur ${\mit\Gamma}$.  Comme il n'y
a, de nouveau, qu'un nombre fini de droites tracées sur $X$, il y a un
ouvert \emph{non vide} (de ${\mit\Gamma}$) de points $x$ pour lesquels
$C(x)$ est une cubique intègre, n'ayant donc pas d'autre singularité
qu'en $x$.  Par ailleurs, comme c'est le cas pour $y$, il y a un
ouvert \emph{non vide} de points $x$ pour lesquels $C(x)$ a deux
directions tangentes distinctes en $x$.  Fixons $x$ dans
l'intersection de ces deux ouverts : $C(x)$ a alors toutes les
propriétés souhaitées, c'est-à-dire qu'il s'agit d'une cubique intègre
ayant pour seule singularité un point double ordinaire en $x$.
\end{proof}

{\footnotesize

En caractéristique $2$, le lemme n'est pas valable : en effet, la
surface cubique (lisse) d'équation $X_0^3 + X_1^3 + X_2^3 + X_3^3 = 0$
est coupée par tout plan tangent soit en une cubique intègre
cuspidale, soit en une droite et une conique tangente à celle-ci, soit
en trois droites concourantes en un point d'Eckardt (elle a $35$
points d'Eckardt).  Il n'existe donc aucune courbe \emph{plane} très
libre sur cette hypersurface cubique (mais il existe, bien sûr, des
courbes très libres tracées dessus, par exemple celle donnée par
$h\colon (U:V) \mapsto (U^3+U^2 V:U^3+U^2 V+V^3:U^2 V+V^3:U V^2)$, qui
a $h^* T_X \cong \mathcal{O}_{\mathbb{P}^1} (2) \oplus
\mathcal{O}_{\mathbb{P}^1} (1)$ puisqu'en contractant six droites bien
choisie elle devient isomorphe à une droite dans le plan éclaté en six
points qui ne passe par aucun de ces points).

Ceci est d'ailleurs le seul exemple possible d'une surface cubique sur
laquelle n'est pas tracée une courbe cubique plane à point double
ordinaire.  En effet, la démonstration du lemme \ref{CubSurfLemma}
vaut encore en caractéristique $2$ dès lors qu'il existe un point à
l'intersection de deux droites qui n'est pas un point d'Eckardt : en
revenant à la démonstration du sous-lemme \ref{SixPointSublemma}, si
tous les points d'intersections sont des points d'Eckardt, en prenant
$(0:0:1) ,\penalty-100 (1:0:1) ,\penalty-100 (1:1:1) ,\penalty-100
(0:1:1)$ pour coordonnées de $P_0, P_1, P_2, P_3$, les points $P_4$ et
$P_5$ doivent être alignés avec les intersections $(1:0:0)$ de $(P_0
P_1)$ et de $(P_2 P_3)$, $(1:1:0)$ de $(P_0 P_2)$ et de $(P_1 P_3)$,
et $(0:1:0)$ de $(P_0 P_3)$ et de $(P_1 P_2)$, donc on peut écrire
$P_4 = (X_4:Y_4:0)$ et $P_5 = (X_5:Y_5:0)$.  Comme $(P_0 P_4)$, $(P_1
P_5)$ et $(P_2 P_3)$ doivent concourir en un point, on doit encore
avoir $Y_4 X_5 + X_4 Y_5 + Y_4 Y_5 = 0$, et, de même, $Y_4 X_5 + X_4
Y_5 + X_4 X_5 = 0$, donc finalement $X_4 X_5 + Y_4 Y_5 = 0$, autrement
dit $P_4 = (1:\zeta:0)$ et $P_5 = (1:\zeta^2:0)$ avec $\zeta$ racine
primitive cubique de l'unité.  Ceci montre que les six points sont
complètement déterminés par l'hypothèse, et il n'y a donc qu'une seule
surface cubique qui vérifie cette propriété.

}

\smallbreak

On peut également facilement trouver des courbes très libres sur des
hypersurfaces cubiques de dimension $\geq 3$ :
\begin{prop}
Soit $X \subseteq \mathbb{P}^n$ une hypersurface cubique lisse de
dimension $n-1 \geq 2$ sur un corps $k$ algébriquement clos (de
caractéristique arbitraire).  Alors il existe $h \colon \mathbb{P}^1
\to X$ très libre sur $X$ (c'est-à-dire $H^1 (\mathbb{P}^1, T_X(-2)) =
0$) et dont l'image est contenue dans un plan si la caractéristique de
$k$ est autre que $2$, et dans un espace de dimension $3$ si elle est
$2$.
\label{CubHyperSurfProp}
\end{prop}
\begin{proof}
On procède par récurrence sur la dimension $n-1$ de l'hypersurface.
Le cas $n-1=2$ constitue le contenu de la
proposition \ref{CubSurfProp} (combiné au lemme \ref{CubSurfLemma}) en
caractéristique différente de deux, et est clair en caractéristique
deux.  Si $n-1\geq 3$, le théorème de Bertini permet de trouver un
hyperplan $H$ tel que l'intersection de $H$ avec $X$ soit une
hypersurface cubique lisse de dimension $n-2$ : l'hypothèse de
récurrence assure alors l'existence de $h \colon \mathbb{P}^1 \to
H\cap X$ très libre dont l'image est contenue dans un plan, et la
proposition \ref{FreenessPropOne} montre que $h$ est libre sur $X$
tout entière.
\end{proof}

\end{document}